\documentclass[11pt]{article}
\newcommand{\version}{August 4, 2008 }

\newcommand{\lanbox}{\hfill \hbox{$\, 
\vrule height 0.25cm width 0.25cm depth 0.01cm
\,$}}

\usepackage{amsthm,amsfonts,amsmath, amscd}
\swapnumbers
\theoremstyle{plain}

\newtheorem{thm}{THEOREM}[section]
\newtheorem{lm}[thm]{LEMMA}

\newtheorem{prop}[thm]{PROPOSITION}
\theoremstyle{definition}
\newtheorem{defi}[thm]{DEFINITION}
\newtheorem{exam}[thm]{EXAMPLE}
\theoremstyle{definition}

\newcommand{\upchi}{\raise1pt\hbox{$\chi$}}
\newcommand{\R}{{\mathord{\mathbb R}}}

\newcommand{\C}{{\mathord{\mathbb C}}}
\newcommand{\Z}{{\mathord{\mathbb Z}}}

\newcommand{\hn}{{\mathord{\widehat{n}}}}

\newcommand{\J}{{\mathcal{J}}}
\newcommand{\tr}{{\rm Tr}}
\renewcommand{\|}{{\Vert}}
\numberwithin{equation}{section}
\begin{document}

\markboth{\scriptsize{EACEHL}}{\scriptsize{EACEHL}}

\title{BRASCAMP-LIEB INEQUALITIES FOR NON-COMMUTATIVE INTEGRATION}
\author{\vspace{5pt} Eric A. Carlen$^1$ and
Elliott H. Lieb$^{2}$ \\
\vspace{5pt}\small{$1.$ Department of Mathematics, Hill Center,}\\[-6pt]
\small{Rutgers University,
110 Frelinghuysen Road
Piscataway NJ 08854-8019 USA}\\
\vspace{5pt}\small{$2.$ Departments of Mathematics and Physics, Jadwin
Hall,} \\[-6pt]
\small{Princeton University, P.~O.~Box 708, Princeton, NJ
  08544}\\
 }
\date{\version}
\maketitle
\footnotetext                                                                         
[1]{Work partially
supported by U.S. National Science Foundation
grant DMS 06-00037.    }                                                          
\footnotetext
[2]{Work partially
supported by U.S. National Science Foundation
grant PHY 06-52854.\\
\copyright\, 2008 by the authors. This paper may be reproduced, in its
entirety, for non-commercial purposes.}

\def\mn{{\bf M}_n}
\def\hn{{\bf H}_n}
\def\hnp{{\bf H}_n^+}
\def\hmnp{{\bf H}_{mn}^+}
\def\h{{\cal H}}
\def\A{{\mathfrak A}}
\def\B{{\mathfrak B}}
\def\CC{{\mathfrak C}}
\def\dd{{\rm d}}
\def\tr{{\rm Tr}}

\begin{abstract}

We formulate a non-commutative analog of the Brascamp-Lieb inequality, and prove it in several
concrete settings.

\end{abstract}

\medskip
\leftline{\footnotesize{\qquad Mathematics subject classification numbers:  47C15, 15A45, 26D15}}
\leftline{\footnotesize{\qquad Key Words: inequalities, traces, non-commutative integration}}

\section{Introduction} \label{intro}

\subsection{Young's inequality in
the context of ordinary Lebesgue integration}

In this paper,  we shall extend the class of  generalized Young's
inequalities known as Brascamp-Lieb inequalities (B-L inequalities)
to an operator algebra setting entailing  non-commutative integration.
\if false We begin by specifying precisely what we mean by a generalized
Young's inequality in the context of ordinary Lebesgue integration.  \fi

The original Young's inequality \cite{Y12} states that for non
negative measurable functions $f_1$, $f_2$ and $f_3$ on $\R$, and
$1\le p_1,p_2,p_3 \le \infty$, with $1/p_1 + 1/p_2 + 1/p_3 = 2$,
\begin{equation}\label{basicy}
\int_{\R^2}f_1(x)f_2(x-y)f_3(y)\dd x\dd y \le \left(\int_\R f_1^{p_1}(t)\dd t\right)^{1/p_1}
\left(\int_\R f_2^{p_2}(t)\dd t\right)^{1/p_2}\left(\int_\R f_3^{p_3}(t)\dd t\right)^{1/p_3}\ .
\end{equation}
Thus, it provides an estimate of the integral of a product of
functions in terms of a product of $L^p$ norms of these functions. The
crucial difference with a H\"older type inequality is that the
integrals on the right are integrals over only $\R$, while the
integrals on the left are integrals over $\R^2$, and none of the three
factors in the product on the left -- $f(x)$, $g(x-y)$ or $h(y)$ --
are integrable (to any power) on $\R^2$.

To frame the inequality in terms that are more amenable to the
generalizations considered here, define the maps $\phi_j: \R^2 \to
\R$, $j=1,2,3$, by
$$\phi_1(x,y) = x\qquad \phi_2(x,y) = x-y\qquad{\rm and}\qquad \phi_3(x,y) = y\ .$$
Then (\ref{basicy}) can be  rewritten as
\begin{equation}\label{basicy2}
\int_{\R^2}\left(\prod_{j=1}^3 f_j\circ \phi_j\right) \dd^2 x \le \prod_{j=1}^3 \left(\int_\R f_j^{p_j}(t)\dd t\right)^{1/p_j}\ .
\end{equation}

There is now no particular reason to limit ourselves to products of only three functions, or to integrals over $\R^2$ and $\R$, or even any Euclidean space for that matter:

\begin{defi} Given measure spaces $(\Omega, \mathcal S , \mu)$ and
  $(M_j, \mathcal M_j , \nu_j)$, $j=1,\dots,N$, not necessarily
  distinct, together with measurable functions $\phi_j:\Omega \to M_j$
  and numbers $p_1,\dots,p_N$ with $1\le p_j \le \infty$, $1 \le j\le
  N$, we say that {\em a B-L inequality holds for
    $\{\phi_1,\dots,\phi_N\}$ and $\{p_1,\dots,p_N\}$} in case there
  is a finite constant $C$ such that
\begin{equation}\label{gy}
\int_\Omega \prod_{j=1}^N f_j\circ \phi_j {\rm d}\mu \le C \prod_{j=1}^N\|f_j\|_{L^{p_j}(\nu_j)}\ 
\end{equation}
\end{defi}
\noindent{holds} whenever each $f_j$ is non negative and measurable on $M_j$, $j=1,\dots,N$.

There are by now many examples. One of the oldest is the original {\em
  discrete Young's inequality}.  In the current notation, this
concerns the case in which $\Omega = \Z^2$ equipped with counting
measure, $N =3$, and each $M_j$ is $\Z$, equipped with counting
measure. Then with
$$\phi_1(m,n) = m\qquad \phi_2(m,n) = m-n\qquad{\rm and} \qquad\phi_3(m,n) = n\ ,$$
(\ref{basicy2}) holds for any three non-negative functions $f_j:\Z \to
\R_+$ under the same conditions on the $p_j$ as in the continuous
case; i.e., $1/p_1 + 1/p_2 + 1/p_3 = 2$.  There is a significant
difference: In the discrete case, the constant $C=1$ is sharp, and
there is equality if and only if one of the $f_j$ is identically zero,
or else $f_1$ vanishes except at some $m_0$, $f_3$ vanishes except at
some $n_0$, and $f_2$ vanishes except at $m_0-n_0$. The inequality
itself is due to Young \cite{Y12}, while the statement about cases of
equality is proved in \cite{HLP}, where the authors also consider
extensions to more than three functions.

In the continuous case, a much wider generalization to more than three
functions was made by  B-L in \cite{BL}, where the sharp constant in Young's
inequality -- which is strictly less than $1$ unless $p_1 = p_2 =1$ --
was obtained, with a proof that the only non-negative functions yielding
equality are certain Gaussian functions.  (This best constant was also
obtained at the same time by Beckner \cite{Beckner}, for three functions.)

These inequalities generalize from $\R$ to $\R^n$.  The complete
generalization to the case in which the $M_j$ are all Euclidean
spaces, but of different dimension, and the $\phi_j$ are linear
transformations from $\R^n$ to $M_j$, was proved by Lieb \cite{L90}.
Again, the maximizers are Gaussians.  Another proof of this
generalized version, together with a reverse form, was obtained by
Barthe \cite{B}, who also provided a detailed analysis of the cases of
equality in the original B-L inequality from \cite{BL}. The cases of
equality in the higher dimensional generalization from \cite{L90} were
analyzed in detail in \cite{BCCT1,BCCT2}.

Examples in which $\Omega$ is the sphere $S^{N-1}$ or the permutation
group ${\cal S}^N$ were proved in \cite{CLL1,CLL2}, and the above
definition of B-L inequalities in arbitrary measure spaces is taken
from \cite{CC}, where a duality between B-L inequalities and
subadditivity of entropy inequalities is proved.

\subsection{A generalized Young's inequality in
the context of non commutative integration}

In non commutative integration theory, as developed by Irving Segal
\cite{S53,S56,S65}, the basic framework is a triple $(\h,\A,\lambda)$
where $\h$ is a Hilbert space, $\A$ is a $W^*$ algebra (a von Neumann
algebra) of operators on $\h$, and $\lambda$ is a positive linear
functional on the finite rank operators in $\A$.  In Segal's picture,
the algebra $\A$ corresponds to the algebra of bounded measurable
functions, and applying the linear positive linear functional
$\lambda$ to a positive operator corresponds to taking the integral of
a positive function. That is,
$$A \mapsto \lambda(A)\qquad{\rm  corresponds\  to}\qquad 
f \mapsto \int_M f\dd \nu\ .$$

Such a triple $(\h,\A,\lambda)$ is called a {\em non commutative
  integration space}.  Certain natural regularity properties must be
imposed on $\lambda$ if one is to get a well-behaved non-commutative
integration theory, but we shall not go into these here as the
examples that we consider are all based on the case in which $\lambda$
is the {\em trace} on operators on $\h$, or some closely related
functional, for which discussion of these extra conditions would be a
digression.

In this operator algebra setting, there are natural non-commutative
analogs of the usual $L^p$ spaces: If $A$ is a finite rank operator in
$\A$, and $1\le q < \infty$, define
$$\|A\|_{q,\lambda} = \left(\lambda(A^*A)^{q/2}\right)^{1/q}\ .$$
This defines a norm (under appropriate conditions on $\lambda$ that
are obvious for the trace), and the completion of the space of finite
rank operator in $\A$ under this norm defines a non-commutative $L^p$
space.  (The completion may contain unbounded operators ``affiliated''
to $\A$.) For more on the general theory of non-commutative
integration, see the early papers \cite{Di53,S53,S65,St59} and the
more recent work in \cite{FK,H,K2,N74}.

To frame an analog of (\ref{gy}) in an operator algebra setting, we
replace the measure spaces by non commutative integration spaces:
$$
(\Omega, \mathcal S , \mu) \longrightarrow (\h, \A,\lambda) \qquad{\rm
  and}\qquad (M_j, \mathcal M_j , \nu_j) \longrightarrow (\h_j,\A_j,
\lambda_j)\qquad j =1,\dots,N \ .$$ 
The right hand side of (\ref{gy})
has an obvious generalization to the operator algebra setting in terms
of the non-commutative $L_p$ norms introduced above.

As for the left hand side of (\ref{gy}), regard $f_j \mapsto f_j\circ
\phi_j$ as a $W^*$ algebra homomorphism (which, restricted to the
$W^*$ algebra $L^\infty(M_j)$, it is), and suppose we are given $W^*$
homomorphisms
$$\phi_j : \A_j \to \A\ .$$
Then each $\phi_j(A_j)$ belongs to $\A$, however in the
non-commutative case, the product of the $\phi_j(A_j)$ depends on
their order in the product, and need not be self-adjoint even -- let
alone positive -- even if each of the $A_j$ are positive.

Therefore, let us return to the left side of (\ref{gy}), and suppose
that each $f_j$ is strictly positive. Then defining
$$h_j = \ln(f_j) \qquad{\rm so \ that}\qquad  f_j\circ \phi_j = e^{h\circ \phi_j}\ ,$$
we can then rewrite (\ref{gy}) as
\begin{equation}\label{gy2}
\int_\Omega \exp\left( \sum_{j=1}^N h_j\circ \phi_j\right) {\rm d}\mu
 \le C \prod_{j=1}^N\| e^{h_j}\|_{L^{p_j}(\nu_j)}\ ,
\end{equation}

We can now formulate our operator algebra analog of  (\ref{gy}):

\begin{defi}\label{ncgy}
Given non commutative integration spaces $(\h,\A,\lambda)$ and  $(\h_j,\A_j,\lambda_j)$, $j=1,\dots,N$,
together with $W^*$ algebra homomorphisms $\phi_j:\A_j\to \A$, $j=1,\dots,N$, and indices
$1\le p_j\le \infty$, $j=1,\dots,N$, a {\em non-commutative B-L inequality holds 
for $\{\phi_1,\dots,\phi_N\}$  and  $\{p_1,\dots,p_N\}$} if there is a finite constant $C$
so that
\begin{equation}\label{intprod3}
\boxed{ \ \lambda\left(\exp\left[ \sum_{j=1}^N \phi_j(H_j)\right]\right) \le C 
\prod_{j=1}^N\left(\lambda_j\exp\left[ p_j H_j \right]\right)^{1/p_j} \ }
\end{equation}
whenever $H_j$ is self-adjoint in $\A_j$, $j = 1,\dots,N$.
\end{defi}

In this paper, we are concerned with determining the indices and the
best constant $C$ for which such an inequality holds, and shall focus
on two examples: The first concerns {\em operators on tensor products
  of Hilbert spaces}, and the second concerns {\em Clifford algebras}.

\subsection{A generalized Young's inequality for tensor products}

\begin{exam}\label{tp}
Let $\h_j$, $j=1,\dots,n$ be separable Hilbert spaces, and let 
Let ${\cal K}$
denote the tensor product
$${\cal K} = \h_1\otimes \cdots \otimes \h_n\ .$$

Define $\A$ to be $\B({\cal K})$, the algebra of bounded operators on
${\cal K}$, and define $\lambda$ to be $\tr$, the trace $\tr$ on
${\cal K}$, so that $(\h,\A,\lambda) = ({\cal K},\B({\cal K}),\tr)$.

For any non empty subset $J$ of $\{1,\dots,n\}$, let ${\cal K}_J$
denote the tensor product
$${\cal K}_J =  \bigotimes_{j\in J}\h_j\ .$$ Define 
$\A_J$ to be $\B({\cal K}_J)$, the algebra of bounded operators on ${\cal K}_J$, and define
$\lambda_J$ be $\tr_J$, the trace on ${\cal K}_J$, so that $(\h_J,\A_J,\lambda_J) = ({\cal K}_J,\B({\cal K}_J),\tr_J)$.

There are natural homomorphisms $\phi_J$ embedding the $2^n-1$
algebras $\A_J$ into $\A$. For instance, if $J = \{1,2\}$,
\begin{equation}\label{embed1}
\phi_{\{1,2\}}(A_1\otimes A_2) = A_1\otimes  A_2\otimes  I_{\h_3}\otimes \cdots  \otimes I_{\h_N}\ ,
\end{equation}
and is extended linearly. 

It is obvious that in case $J\cap K = \emptyset$ and $J\cup K
=\{1,\dots,n\}$, then for all $H_J \in \A_J$ and $H_K \in \A_K$,
\begin{equation}\label{gt1}
\tr\left(e^{H_J+H_K}\right) = \tr_J\left(e^{H_J}\right)\tr_K\left(e^{H_K}\right)\ ,
\end{equation}
but things are more interesting when $J\cap K\ne \emptyset$ and $J$
and $K$ are both proper subsets of  $\{1,\dots,n\}$.  If $H_J$ and
$H_K$ do not commute, which is the typical situation for $J\cap K\ne
\emptyset$, one can estimate  the left hand side of (\ref{gt1}) by first
applying the Golden--Thompson inequality \cite{Go,T}, which says that
for self-adjoint operators $H_J$ and $H_K$, $$\tr\left(e^{H_J+H_K}\right)
\le \tr\left(e^{H_J}e^{H_K}\right) \ .$$ One might then apply  H\"older's
inequality -- but  if $J$ and $K$ are proper subsets of $\{1,\dots,n\}$,
this will yield a finite bound if and only if all of the Hilbert
spaces whose indices are not included in {\em both} $J$ and $K$ are
finite dimensional. Even then, the bound depends on the dimension in an
unpleasant way.  The non-commutative B-L Inequalities provided by the
next theorem  do not have this defect.

\begin{thm}\label{simpletp} Let $J_1,\dots, J_N$ be $N$ non empty
  subsets of $\{1,\dots,n\}$ For each $i \in \{1,\dots,n\}$, let
  $p(i)$ denote the number of the sets $J_1,\dots,J_N$ that contain
  $i$, and let $p$ denote the minimum of the $p(i)$. Then, for
  self-adjoint operators $H_j$ on ${\cal K}_{J_j}$, $j=1,\dots,N$,
\begin{equation}\label{tpgy}
\boxed{\ \tr \left(\exp\left[ \sum_{j=1}^N \phi_{J_j}(H_j)\right]\right) \le  
\prod_{j=1}^N\left(\tr_{J_j}\  e^{ q H_j }\right)^{1/q} \ }
\end{equation}
for $q=p$ (and hence all $1 \le q \le p$), while for all $q>p$, it is
possible for the left hand side to be infinite, while the right hand
side is finite.
\end{thm}

Note that in Theorem~\ref{simpletp}, the constant $C$ in Definition
(\ref{ncgy}) is $1$.  The fact that the constant $C=1$ is best
possible, and that the inequality cannot hold for $q>p =
\min\{p(1),\dots,p(N)\}$ is easy to see by considering the case that
each $\h_j$ has finite dimension $d_j$, and $H_j=0$ for each $j$. Then
$$\tr \left(\exp\left[ \sum_{j=1}^N \phi_{J_j}(H_j)\right]\right) = \prod_{j=1}^n d_j \qquad{\rm and}\qquad
\prod_{j=1}^N \left(\tr_{J_j} e^{ q H_j }\right)^{1/q}  = \prod_{j=1}^N \prod_{k\in J_j} d_k^{1/q} =
\prod_{j=1}^n  d_j^{p(j)/q}\ .$$
We will prove the inequality (\ref{tpgy}) for $q=p$  in Section~\ref{SSA}.

As an example, consider the case in which $n=6$, $N =3$ and
$$J_1 = \{1,2,3\}\quad J_2=\{3,4,5\} \qquad{\rm and}\quad J_3 = \{5,6,1\}\ .$$
Here,  $p=1$, and hence
\begin{equation}\label{intprod4c}
\tr \left(\exp\left[ \sum_{j=1}^3 \phi_{J_j}(H_j)\right]\right) \le  
\prod_{j=1}^3\left(\tr_{J_j}\  e^{  H_j }\right) \ .
\end{equation}
The inequality (\ref{intprod4c}) can obviously be extended to larger
tensor products, and has an interesting statistical mechanical
interpretation as a bound on the {\em partition function} of a
collection of interacting spins in terms of a product of partition
functions of simple constituent sub-systems.

To estimate the left side of (\ref{intprod4c}) without using
Theorem~\ref{simpletp}, one might use the Golden-Thompson inequality
and then Schwarz's inequality to write
$$\tr \left(\exp\left[ \sum_{j=1}^3 \phi_{J_j}(H_j)\right]\right) \le 
\tr \left( e^{ \phi_1(H_1)+ \phi_3(H_3)} e^{\phi_2(H_2)} \right) \le
\left(\tr\ e^{2[\phi_1(H_1)+ \phi_3(H_3)]}\right)^{1/2}\left(\tr\
  e^{2\phi_2(H_2)}\right)^{1/2}\ .$$ While the $L^2$ norms are an
improvement over the $L^1$ norms in (\ref{intprod4c}), the traces are
now over the entire tensor product space. Thus, for example,
$$\left(\tr\ e^{2\phi_2(H_2)}\right)^{1/2} = (d_1d_2d_6)^{1/2}\left(\tr_{J_2}\ e^{2 H_2 }\right)^{1/2}$$
where $d_j$ is the dimension of Hilbert space $\h_j$.   This dimension dependence may be unfavorable  if any of the dimensions is large.

\end{exam}

\subsection{A generalized Young's inequality in Clifford algebras}

Our next example concerns Clifford algebras, which as Segal emphasized
\cite{S56}, allow one to represent Fermion Fock space as an $L^2$
space -- albeit a non-commutative $L^2$ space, but still with many of
the advantages of having a Hilbert space represented as a function
space, as in the usual Schr\"odinger representation in quantum
mechanics.

In the finite dimensional setting, with $n$ degrees of freedom, one
starts with $n$ operators $Q_1,\dots,Q_n$ on some Hilbert space $\h$
that satisfy the {\em canonical anticommutation relations}
$$Q_iQ_j + Q_jQ_i = 2\delta_{i,j} I\ .$$
One can concretely construct such operators acting on $\h =
(\C^2)^{\otimes n}$, the $n$--fold tensor product of $\C^2$ with
itself; see \cite{BW}.  The Clifford algebra $\CC$ is the operator
algebra on $\h$ that is generated by $Q_1,\dots,Q_n$.

The Clifford algebra $\CC$ itself is $2^n$ dimensional. In fact, let
$\alpha = (\alpha_1,\dots,\alpha_n)$ be a {\em Fermionic multi-index},
which means that each $\alpha_j$ is either $0$ or $1$. Then define
\begin{equation}\label{multii}
Q^\alpha = Q_1^{\alpha_1}Q_2^{\alpha_2}\cdots Q_n^{\alpha_n}\ .
\end{equation}
it is easy to see that the $2^n$ operators $Q^\alpha$ are a basis for
the Clifford algebra, so that any operator $A$ in $\CC$ has a unique
expression
$$A = \sum_{\alpha}x_\alpha Q^\alpha\ .$$
The linear functional $\tau$  on $\CC$ is defined by
\begin{equation}\label{taudef}
\tau\left(\sum_{\alpha}x_\alpha Q^\alpha\right) = x_{(0,\dots,0)}\ .
\end{equation}
That is, $\tau$ acting on $A$ picks off the coefficient of the
identity in $A= \sum_{\alpha} x_\alpha Q^\alpha$.  It turns out that
when the Clifford algebra is constructed in the way described here, as
an algebra operators on the $2^n$ dimensional space $\h$, $\tau$ is
nothing other than the normalized trace:
$$\tau(A) = \frac{1}{2^n}\tr_\h(A)\ .$$
Hence $\tau$ is a positive linear functional, and $((\C^2)^{\otimes
  n},\CC, \tau)$ is a non commutative integration space in the sense
of Segal.

Clifford algebras have infinitely many subalgebras that are also
Clifford algebras of lower dimension. This is in contrast to the
setting described in Example~\ref{tp}, in which the only natural
subalgebras are the $2^n-1$ subalgebras corresponding to the $2^n-1 $
non empty subsets of the index set $\{1,\dots,n\}$.

To describe these subalgebras, let ${\cal J}$ be the {\em canonical
  injection} of $\R^n$ into $\CC$, which is given by
\begin{equation}\label{caninj}{\cal J}((x_1,\dots,x_n)) = \sum_{j=1}^n x_jQ_j\ .
\end{equation}
If $x$ and $y$ are any two vectors in $\R^n$, it is easy to see from
the canonical anticommutation relations that
$$ ({\cal J}(x))({\cal J}(y)) = 2(x\cdot y) I\ .$$
Hence if $V$ is {\em any} $m$ dimensional subspace of $\R^n$, and
$\{u_1,\dots,u_m\}$ is {\em any} orthonormal basis for $V$, the $m$
operators
$${\cal J}(u_1), \dots, {\cal J}(u_m)$$
again satisfy the canonical anticommutation relations, and generate a
subalgebra of $\CC$ that is denoted by $\CC(V)$, and referred to as
{\em the Clifford algebra over $V$}. In the same vein, it is
convenient to refer to $\CC$ itself as the {\em Clifford algebra over
  $\R^n$}.  Obviously, $\CC(V)$ is naturally isomorphic to
$\CC(\R^m)$, and for $A \in \CC(V)$ one may compute $\tau(A)$ using
either the normalized trace $\tau$ inherited from $\CC$, or the
normalized trace $\tau_V$ induced by the identification of $\CC(V)$
with $\CC(\R^m)$.

As Segal emphasized, $((\C^2)^{\otimes n},\CC,\tau)$ is in many way a
non-commutative analog of the Gaussian measure space
$(\R^n,\gamma(x)\dd x)$ where
\begin{equation}\label{gauss}
\gamma(x) = \frac{1}{(2\pi)^{n/2}}e^{-|x|^2/2}\ .
\end{equation}
In fact, just as orthogonality implies independence in
$(\R^n,\gamma(x)\dd x)$, if $V$ and $W$ are two orthogonal subspaces
of $\R^n$, and if $A\in \CC(V)$ and $B\in \CC(W)$, then
$$\tau(AB) = \tau(A)\tau(B)\ .$$
The results we prove here reenforce this analogy.  We are now ready to
introduce our next example:

\begin{exam}\label{cliff}
  For some $n>1$, let $\A$ be the Clifford algebra over $\R^n$ with
  its usual inner product, and let $\A$ be equipped with its unique
  tracial state $\tau$, which is simply the normalized trace.

  For each $j=1,\dots,N$, let $V_j$ be a subspace on $\R^n$, and let
  $\A_j$ be $\CC(V_j)$, the Clifford algebra over $V_j$ with the inner
  product $V_j$ inherits from $\R^n$. Let $\A_j$ be equipped with its
  unique tracial state $\tau_j$.  The natural embedding of $V_j$ into
  $\R^n$ induces a homomorphism of $\A_j$ into $\A$, and we define
  this to be $\phi_j$.  In this setting, we shall prove

  \begin{thm}\label{cy} Let $V_1,\dots,V_N$ be $N$ subspaces of
    $\R^n$, and let $\A_j$ be the Clifford algebra over $V_j$ with the
    inner product $V_j$ inherits from $\R^n$, and let $\A_j$ be
    equipped with its unique tracial state $\tau_j$.  Let $\phi_j$ be
    the natural homomorphism of $\A_j$ into $\A$ induced by the
    natural embedding of $V_j$ into $\R^n$.  Then
\begin{equation}\label{intprod5b}
\boxed{\ 
\tau \left(\exp\left[ \sum_{j=1}^N \phi_j(H_j)\right]\right) \le  
\prod_{j=1}^N\left(\tau_j\  e^{ p_j H_j }\right)^{1/p_j} \ }
\end{equation}
for all self-adjoint operators $H_j\in \A_j$ if and only if
\begin{equation}\label{fjc}
\sum_{j=1}^N\frac{1}{p_j}P_j \le I_{\R^n}\ .
\end{equation}
where 
 $P_j$ is the orthogonal projection onto
$V_j$ in $\R^n$. 
\end{thm}

In the special case in which ${\rm dim}(V_j) = 1$ for each $j$,
(\ref{intprod5b}) reduces to an interesting inequality for the
hyperbolic cosine. Indeed, let $u_j$ be one of the two unit vectors in
$V_j$.

Then, with $u_j\otimes u_j$ denoting the orthogonal projection onto the span of $u_j$, (\ref{fjc})
reduces to
\begin{equation}\label{fjc1d}
\sum_{j=1}^N\frac{1}{p_j}u_j\otimes u_j  \le I_{\R^n}\ .
\end{equation}
The greater simplification, however, is that in this case, the space
of self-adjoint operators in each $\A_j$ is just two dimensional, and
with ${\cal J}$ denoting the canonical injection defined in
(\ref{caninj}),
$$H_j = a_jI + b_j{\cal J}(u_j)$$
for some real numbers $a_j$ and $b_j$.  Then
$$\sum_{j=1}^N H_j = \left(\sum_{j=1}^N a_j\right)I  + {\cal J}\left(\sum_{j=1}^N b_ju_j\right)\ .$$
This operator has exactly two eigenvalues,
$$\left(\sum_{j=1}^N a_j\right) \pm  \left|\sum_{j=1}^N b_ju_j\right|\ $$
with equal multiplicities.  

Likewise, $p_jH_j$ has exactly two eigenvalues $p_ja_j \pm p_jb_j$
with equal multiplicities.  Hence, in this case, (\ref{intprod5b})
reduces to
\begin{equation}\label{cosh1}
\cosh\left(\left|{\textstyle \sum_{j=1}^N b_ju_j}\right|\right) \le \prod_{j=1}^N\left(\cosh(p_jb_j)\right)^{1/p_j}\quad{\rm for\ all}\quad (b_1,\dots,b_N) \in \R^N\ ,
\end{equation}
which, according to the theorem, must hold whenever (\ref{fjc1d}) is
satisfied.  (The $a_j$'s make the same contribution to both sides, and
may be cancelled away.)  Taking the logarithm of both sides, this can
be rewritten as
\begin{equation}\label{cosh2}
\ln\cosh\left(\left|{\textstyle \sum_{j=1}^N b_ju_j}\right|\right) \le \sum_{j=1}^N\frac{1}{p_j}
\ln\cosh(p_jb_j) \quad{\rm for\ all}\quad (b_1,\dots,b_N) \in \R^N\ ,
\end{equation}
and this inequality must hold whenever the unit vectors
$\{u_1,\dots,u_N\}$ and the positive numbers $\{p_1,\dots,p_N\}$ are
such that (\ref{fjc1d}) is satisfied.

Later on, we shall give an elementary proof of this inequality, and
hence an elementary proof of Theorem~\ref{cy} when each $V_j$ is one
dimensional. Our proof of the other cases is less than elementary, and
even our elementary proof of (\ref{cosh2}) is less than direct.

\end{exam}

\section{Subadditivty of Entropy and Generalized Young's Inequalities}\label{subadd}

In the examples we have introduced in the previous section, the
positive linear functionals $\lambda$ under consideration are either
traces or normalized traces.  Throughout this section, we assume that
our non commutative integration spaces $(\h,\A,\lambda)$ are based on
{\em tracial} positive linear functionals $\lambda$. That is, we
require that for all $A,B \in \A$,
 $$\lambda(AB) = \lambda(BA)\ .$$

 In such a non commutative integration space $(\h,\A,\lambda)$, a {\em
   probability density} is a non negative element $\rho$ of $\A$ such
 that $\lambda(\rho) =1$. Indeed, the tracial property of $\lambda$
 ensures that
$$\lambda(\rho A) = \lambda(A\rho) = \lambda(\rho^{1/2}A\rho^{1/2})$$
so that $A \mapsto \lambda(\rho A)$ is a positive linear functional that is $1$ on the identity.

Now suppose we have $N$ non-commutative integration spaces
$(\h_j,\A_j,\lambda_j)$ and $W^*$ homomorphism $\phi_j: \A_j \to \A$.
Then these homomorphisms induce maps from the space of probability
densities on $\A$ to the spaces of probability densities on the
$\A_j$:
 
For any probability density $\rho$ on $(\A, \lambda)$, let $\rho_j$ be
the probability density on $(\A_j, \lambda_j)$ by
$$\lambda_j(\rho_j A) = \lambda(\rho \phi_j(A))$$
for all $A \in \A_j$. 

For example, in the setting of Example~\ref{tp}, $\rho_{J_j}$ is just
the partial trace of $\rho$ over $\bigotimes_{k\in J_j^c}\h_{k}$
leaving an operator on $\bigotimes_{k\in J_j}\h_{k}$.  In the Clifford
algebra setting of Example~\ref{cliff}, $\rho_j$ is simply the
orthogonal projection of $\rho$ in $L^2(\CC,\tau)$ onto $\CC(V_j)$,
which is also known as the conditional expectation \cite{Um1} of
$\rho$ given $\CC(V_j)$.

In this section, we are concerned with the relations between the {\em entropies} of $\rho$ and the $\rho_1,\dots,\rho_N$.
The entropy of a probability density $\rho$, $S(\rho)$, is defined by 
$$S(\rho) = -\lambda(\rho \ln \rho)\ .$$
Evidently, the entropy functional is  concave on the set of probability densities.

\begin{defi}\label{gsa}
  Given tracial non-commutative integration spaces $(\h,\A,\lambda)$
  and $(\h_j,\A_j,\lambda_j)$, $j=1,\dots,N$, together with $W^*$
  algebra homomorphisms $\phi_j:\A_j\to \A$, $j=1,\dots,N$, and
  numbers $1\le p_j\le \infty$, $j=1,\dots,N$, a {\em generalized
    subadditivity of entropy inequality} holds if there is a finite
  constant $C$ so that
\begin{equation}\label{genen}
\boxed{\ 
\sum_{j=1}^N \frac{1}{p_j} S(\rho_j) \ge S(\rho)- \ln C\ }
\end{equation}
for all probability densities $\rho$ in $\A$.
\end{defi}

It turns out that for tracial non-commutative integration spaces,
generalized subadditivity of entropy inequalities and B-L inequalities
are dual to one another, just as they are in the commutative case
\cite{CC}, so that if one holds, so does the other, with the same
values of $p_1,\dots,p_N$ and $C$.  The following is in fact a direct
non-commutative analog of the main theorem of \cite{CC}.

\begin{thm}\label{equiv} Let 
  $(\h,\A,\lambda)$ and $(\h_j\A_j,\lambda_j)$, $j=1,\dots,N$, be
  tracial non-commutative integration spaces. Let $\phi_j:\A_j\to \A$,
  $j=1,\dots,N$ be $W^*$ algebra homomorphisms. Then for any numbers
  $1\le p_j\le \infty$, $j=1,\dots,N$, and any finite constant $C$,
  the generalized subadditivity of entropy inequality (\ref{genen}) is
  true for all probability densities $\rho$ on $\A$ if and only if the
  non-commutative B-L inequality (\ref{intprod3}) is true for all
  self-adjoint $H_j\in \A_j$, $j=1,\dots,N$, with the same
  $p_1,\dots,p_N$ and the same $C$.
\end{thm}

As a consequence of Theorem~\ref{equiv}, one strategy for proving a non-commutative 
B-L inequality
is to prove the corresponding generalized subadditivity of entropy  inequality. We shall see in our examples that this is an effective strategy; indeed, this is how we prove Theorems~\ref{simpletp}
and \ref{cy}.

In the current tracial context, the proof of Theorem~\ref{equiv} is a
direct adaptation of the proof of the corresponding result in the
context of Lebesgue integration given in \cite{CC}.  It turns on a
well--known formula for the Legendre transform of the entropy. For
completeness, we give this formula in Lemma~\ref{leg} below. Before
stating the lemma, it is convenient to extend the definition of $S$ to
all of $\A_{\rm sa}$, the subspace of self-adjoint elements of $\A$,
as follows:

\begin{equation}\label{hdef}
S(A) = 
\begin{cases}
-\lambda(A \ln A) & \text{if $A\ge 0$ and $\lambda(A) = 1$,}
\\
-\infty & \text{otherwise.}
\end{cases}
\end{equation}

\begin{lm}\label{leg} Let $\A$ be $\B(\h)$, the algebra of bounded operators on a separable Hilbert space $\h$.
  Let $\lambda$ denote either the trace $\tr$ on $\h$, or, if $\h$ is
  finite dimensional, the normalized trace $\tau$.  Then for all $A\in
  \A_{\rm sa}$,
\begin{equation}\label{forone}
-S(A) = \sup_{H \in \A_{\rm sa}}\left\{ 
\lambda(AH) - \ln\left(\lambda\left(e^H \right)\right)\right\} \ .
\end{equation}
The supremum is an attained maximum if and only if $A$ is a strictly
positive probability density, in which case it is attained at $H$ if
and only if $H = \ln A + cI$ for some $c\in \R$.  Consequently, for
all $H\in \A_{\rm sa}$,
\begin{equation}\label{fortwo}
\ln\left(\lambda\left(e^H\right)\right) = \sup_{A \in \A_{\rm sa}}\left\{ \lambda(AH)+ S(A)\right\} \ .
\end{equation}
The supremum is a maximum at all points of the domain of $
\ln\left(\lambda\left(e^H\right)\right)$, in which case it is attained
only at the single point $A = e^H/(\lambda(e^H))$.
\end{lm}

\noindent{\bf Proof:} We consider first the case that $\lambda = \tr$,
and $\h$ has finite dimension $d$. To see that the supremum is
$\infty$ unless $0 \le A \le I$, let $c$ be any real number, and let
$u$ be any unit vector. Then let $H$ be $c$ times the orthogonal
projection onto $u$.  For this choice of $H$,
$$\lambda(AH) - \ln\left(\lambda\left(e^H \right)\right) = c\langle u,Au\rangle  - \ln (e^c + (d-1))\ .$$
If $\langle u,Au\rangle<0$, this tends to $\infty$ as $c$ tends to $-\infty$. If 
$\langle u,Au\rangle>1$, this tends to $\infty$ as $c$ tends to $\infty$. Hence we need only consider $0 \le A \le I$.
Next, taking $H = cI$, $c\in \R$,
$$\lambda(AH) - \ln\left(\lambda\left(e^H \right)\right) = c\lambda(A)  - c - \ln (d)\ .$$
Unless $\lambda(A) = 1$, this tends to $\infty$ as $c$ tends to $\infty$. Hence we need only consider 
the case that $A$ is a  density matrix $\rho$. 

Let $\rho$ be any density matrix on $\h$ and let $H$ be any
self-adjoint operator such that $\tr(e^H) < \infty$.  Then define the
density matrix $\sigma$ by
$$\sigma = \frac{e^H}{\tr(e^H)}\ .$$
Then, by the positivity of the relative entropy,
$$\tr( \rho\ln \rho -  \rho \ln \sigma) \ge 0$$
with equality if and only if $\sigma = \rho$. But by the definition of
$\sigma$, this reduces to
$$\tr (\rho \ln \rho) \ge \tr(\rho H) - \ln\left(\tr\left(e^H\right)\right)\ ,$$
with equality if and only if $H = \ln \rho$.  From here, there rest is
very simple, including the treatment of the normalized trace.. \lanbox
\medskip

Petz \cite{P88} has shown how to extend Lemma~\ref{leg} to a much more
general context, and his result can be used to extend the validity
Theorem ~\ref{equiv} beyond the tracial case. However, since the
examples in which we prove the generalized subadditivity inequality
here are tracial, we shall not go into this.

\medskip
\noindent{\bf Proof of Theorem~\ref{equiv}:} 
Suppose first that the non-commutative B-L inequality (\ref{intprod3})
holds.  Then, for any probability density $\rho$ in $\A$, and any
self-adjoint $H_j \in \A_j$, $j=1,\dots,N$, apply (\ref{forone}) with
$A = \rho$ and $H = \sum_{j=1}^N\phi_j(H_j)$ to obtain
\begin{eqnarray}
-S(\rho) &\ge&  \lambda\left(\rho\left[\sum_{j=1}^N\phi_j(H_j)\right]\right) - 
\ln\left[\lambda \left( \exp\left[\sum_{j=1}^N\phi_j(H_j)\right]\right)\right]\nonumber\\
&=&\sum_{j=1}^N\lambda_j(\rho_j H_j) - \ln\left[\lambda \left( \exp\left[\sum_{j=1}^N\phi_j(H_j)\right]\right)\right]\nonumber\\
&\ge&\sum_{j=1}^N\lambda_j(\rho_j H_j) - \ln\left[C\prod_{j=1}^N
\lambda_j \left( e^{p_jH_j}\right)^{1/p_j}\right]\nonumber\\
&=&\sum_{j=1}^N\frac{1}{p_j}\left[\lambda_j(\rho_j [p_jH_j]) - \ln\left(
\lambda_j \left( e^{[p_jH_j]}\right)\right)\right] - \ln C\ .\nonumber\\
\end{eqnarray}
The first inequality here is (\ref{forone}), and the second is 
the non-commutative B-L inequality (\ref{intprod3}).

Now choosing $p_jH_j$ to maximize $\lambda_j(\rho_j [p_jH_j]) -
\ln\left( \lambda_j \left( e^{[p_jH_j]}\right)\right)$, we get from
(\ref{forone}) once more that
$$\lambda_j(\rho_j [p_jH_j]) - \ln\left(
\lambda_j \left( e^{[p_jH_j]}\right)\right) = -S(\rho_j) = \lambda_j(\rho_j\ln \rho_j)\ .$$
Thus, we have proved (\ref{genen}) with the same $p_1,\dots,p_N$ and $C$ that we had in 
(\ref{intprod3}).

Next, suppose that (\ref{genen}) is true. We shall show that in this
case, the non-commutative B-L inequality (\ref{intprod3}) holds with
the same $p_1,\dots,p_N$ and $C$. To do this, let the self-sadjoint
operators $H_1,\dots,H_N$ be given, and define
$$\rho  = \left[\lambda \left( \exp\left[\sum_{j=1}^N\phi_j(H_j)\right]\right)\right]^{-1}
\exp\left[\sum_{j=1}^N\phi_j(H_j)\right]\ .$$
Then by Lemma~\ref{leg},
\begin{eqnarray}\label{back}
\ln \left[\lambda \left( \exp\left[\sum_{j=1}^N\phi_j(H_j)\right]\right)\right] &=& 
\lambda\left(\rho \left[\sum_{j=1}^N\phi_j(H_j)\right]\right) + S(\rho)\nonumber\\
&=& \sum_{j=1}^N\lambda_j \left[\rho_j H_j\right] + S(\rho)\nonumber\\
&\le& \sum_{j=1}^N \frac{1}{p_j}\left[\lambda_j \left[\rho_j (p_jH_j)\right]  + S(\rho_j)\right] + \ln C \nonumber\\
&\le& \sum_{j=1}^N \frac{1}{p_j}\ln \left[\lambda_j \left(\exp( p_jH_j)\right) \right] + \ln C \nonumber\\
\end{eqnarray}
The first inequality is the generalized subadditivity of entropy
inequality (\ref{genen}), and the second is (\ref{fortwo}).

Exponentiating both sides of (\ref{back}), we obtain (\ref{intprod3})
with the same $p_1,\dots,p_N$ and $C$ that we had in (\ref{genen}).
\lanbox \medskip


\section{Proof of the generalized subadditivity of entropy inequality for tensor products of Hilbert spaces}\label{SSA}

The crucial tool that we use in this section is the {\em strong
  subadditivity of the entropy} \cite{LR}, which we now recall in a
formulation that is suited to our purposes.

Suppose, as in Example~\ref{tp}, that we are given $n$ separable
Hilbert spaces $\h_1,\dots,\h_n$.  As before, let ${\cal K}$ denote
their tensor product, and for any non empty subset $J$ of
$\{1,\dots,n\}$, let ${\cal K}_J$ denote $\bigotimes_{j\in J}\h_j$.

For a density matrix $\rho$ on ${\cal K}$, and any non empty subset
$J$ of $\{1,\dots,n\}$, define $\rho_J = \tr_{J^c}\rho$ to be the
density matrix on ${\cal K}_J$ induced by the natural injection of
$\B({\cal K}_J)$ into $\B({\cal K})$. As noted above, $\rho_J$ is
nothing other than the partial trace of $\rho$ over the complementary
product of Hilbert spaces, $\bigotimes_{j\notin J}\h_j$.

The strong subadditivity of entropy is expressed by the inequality
stating that for all nonempty $J,K \subset \{1,\dots,n\}$,
\begin{equation}\label{ssa}
S(\rho_{J}) + S(\rho_{K})  \ge S(\rho_{J\cup K}) + S(\rho_{J\cap K})\ .
\end{equation}
In case $\J\cap K = \emptyset$, it reduce to the
ordinary subadditivity of the entropy, which  is the elementary inequality
\begin{equation}\label{rsa}
S(\rho_J) + S(\rho_K) \ge S(\rho_{J\cup K})\qquad{\rm for}\quad J\cap K = \emptyset\ .
\end{equation}

Combining these, we have
\begin{eqnarray}
S(\rho_{\{1,2\}}) + S(\rho_{\{2,3\}})+ S(\rho_{\{3,1\}}) &\ge& S(\rho_{\{1,2,3\}})+ S(\rho_{\{2\}}) + S(\rho_{\{1,3\}})\nonumber\\ 
&\ge &
2S(\rho_{\{1,2,3\}})\ ,\nonumber\\
\end{eqnarray}
where the first inequality is the strong subadditivity (\ref{ssa}) and
the second is the ordinary subadditivity (\ref{rsa}).  Thus, for $n=3$
and $J_1= \{1,2\}$, $J_2= \{2,3\}$ and $J_3= \{3,1\}$, we obtain
$$\frac{1}{2}\sum_{j=1}^N S(\rho_{J_j}) \ge S(\rho)\ .$$
In this example, each index $i \in \{1,1,3\}$ belonged to exactly two
of the set $J_1$, $J_2$ and $J_3$, and this is the source of the facto
of $1/2$ in the inequality. The same procedure leads to the following
result:

\begin{thm}\label{stpen} Let $J_1,\dots, J_N$ be $N$ non empty subsets of $\{1,\dots,n\}$  For each $i \in \{1,\dots,n\}$, let $p(i)$
  denote the number of the sets $J_1,\dots,J_N$ that contain $i$, and
  let $p$ denote the minimum of the $p(i)$. Then
\begin{equation}\label{sentp}
\frac{1}{p}\sum_{j=1}^N S(\rho_{J_j}) \ge S(\rho)
\end{equation}
for all density matrices $\rho$ on ${\cal K} = \h_1\otimes \cdots\otimes \h_n$.
\end{thm}

\noindent{\bf Proof:} Simply use strong subadditivty to combine
overlapping sets to produce as many ``complete'' sets as possible, as
in the example above. Clearly, there can be no more than $p$ of these.
If $p(i) > p$ for some indices $i$, there will be ``left over''
partial sets. The entropy is always non negative, and therefore,
discarding the corresponding entropies gives us $\sum_{j=1}^N
S(\rho_{J_j}) \ge pS(\rho)$, and hence the inequality. \lanbox

It is now a very simple matter to prove Theorem~\ref{simpletp}:

\noindent{\bf Proof of Theorem~\ref{simpletp}:} By the remarks made
after the statement of the theorem, all that remains to be proved is
the inequality (\ref{tpgy}) for $q=p$. However, this follows directly
from Theorem~\ref{equiv} and Theorem~\ref{stpen}. \lanbox


\section{On the generalized Young's  inequality with a Gaussian reference measure}\label{gaussY}

Before turning to the proof of our non-commutative B-L inequality in
Clifford algebras, we discuss the commutative case in which the
reference measures is Gaussian. We do this here for two reasons:
First, as noted, a Clifford algebra $\CC$ with its normalized trace
$\tau$ is a non commutative analog of a Gaussian measure space. This
analogy is strong enough that we shall be able to pattern our analysis
in the Clifford algebra case on an analysis of the Gaussian case.

Second, the Gaussian inequality is of interest in itself, and seems
not to have been fully studied before.  Suppose that $V_1,\dots,V_N$
are $N$ non zero subspaces of $\R^n$, and for each $j$, define $\phi_j
= P_j$ to be the orthogonal projection of $\R^n$ onto $V_j$.  Equip
$\R^n$ and equip each $V_j$ with Lebesgue measure. Then the problem of
determining for which sets of indices $\{p_1,\dots,p_N\}$ there exists
a finite constant $C$ so that (\ref{gy}) holds for all non-negative
measurable functions $f_j$ on $V_j$, $j=1,\dots,N$ is highly non
trivial, and has only recently been fully solved \cite{BCCT1,BCCT2}.
Moreover, determining the value of the best constant $C$ for those
choices of $\{p_1,\dots,p_N\}$ is still a challenging finite
dimensional variational problem for which there is no general explicit
solution.

In contrast, suppose we are given a non-degenerate Gaussian measure on
$\R^n$. It will be convenient to take the covariance matrix of the
Gaussian to define the inner product, so that the Gaussian becomes a
unit Gaussian.  For each positive integer $m$, define $\gamma_m(x) =
(2\pi)^{-m/2} e^{-|x|^2/2}$ on $\R^m$. Then equipping $\R^n$ with the
measure $\gamma_n(x)\dd x$ and equipping each $V_j$ with the
$\gamma_{d_j}(x)\dd x$, $d_j$ being the dimension of $V_j$, it turns
out that there is a {\em very simple necessary and sufficient
  condition} on the indices $\{p_1,\dots,p_N\}$ for the constant $C$
to be finite, and better yet, {\em the best constant $C$ is always $1$
  whenever it is finite}:

\begin{thm}\label{gaussbl}  Let $V_1,\dots,V_N$ be $N$ non zero subspaces of $\R^n$, 
  and for each $j$, and let $d_j$ denote the dimension of $V_j$.
  Define $P_j$ to be the orthogonal projection of $\R^n$ onto $V_j$.
  Given the numbers $p_j$, $1\le p_j < \infty$ for $j = 1,\dots,N$,
  there exists a finite constant $C$ such that
\begin{equation}\label{blg}
\int_{\R^n} \prod_{j=1}^N f_j\circ P_j(x) \gamma_n(x){\rm d}x \le C
\prod_{j=1}^N\left(\int_{V_j}f_j^{p_j}(y)\gamma_{d_j}(y){\rm d}y\right)^{1/p_j}\ 
\end{equation}
holds for all non-negative $f_j$ on $V_j$, $j=1,\dots,N$,  if and only if
\begin{equation}\label{fja}
\sum_{j=1}^N \frac{1}{p_j} P_j \leq {\rm Id}_{\R^n}
\end{equation}
and in this case, $C=1$.
\end{thm}

\medskip

We hasten to point out that this theorem is partially known. In the
special case that each of the subspaces $V_j$ is one dimensional,
Barthe and Cordero-Erausquin \cite{BC}, have the {\em sufficiency} of
the condition (\ref{fja}) which reduces to
\begin{equation}\label{geo}\sum_{j=1}^N \frac{1}{p_j}u_j\otimes u_j = {\rm Id}_{\R^n}
\end{equation}
with each $u_j$ being a unit vector spanning $V_j$. They did this as
an intermediate step in a short proof of the {\em Lebesgue measure}
version of the B-L inequality under the condition (\ref{geo}) -- the
so-called {\em geometric case}. Perhaps because their main focus was
the Lebesgue measure case, in which (\ref{geo}) is not a necessary
condition for finiteness of the constant $C$, they did not address the
necessity of this condition in the Gaussian case.

Indeed, the inequality (\ref{blg}) is equivalent to its Lebesgue
measure analog, which is known to hold with the constant $C=1$ under
the condition (\ref{fja}) \cite{BCCT1,BCCT2}. To see this, define
$g_1,\dots,g_N$ by
$$g_j(y) = f_j(y)(\gamma_j(y))^{1/d_j}\qquad j = 1,\dots,N\ .$$
As noted in \cite{BC}, this change of variable allows one to pass back
and forth between the Gaussian and Lebesgue measure version of the B-L
inequality -- under the condition (\ref{fjc}).

Nonetheless, it is worthwhile to give a proof of Theorem~\ref{gaussbl}
here for two reasons: First, it may be surprising that the condition
(\ref{fjc}) is necessary for the inequality to hold with any finite
constant at all.  Second, the proof we will give of sufficiency of the
condition (\ref{fjc}) serves as a model for the proof of the
corresponding theorem in the Clifford algebra case that we consider in
the next section.

In proving Theorem~\ref{gaussbl}, our first step is to pass to the
problem of proving a generalized subadditivity inequality.  Because
the commutative version of Theorem~\ref{equiv} has been proved in
\cite{CC}, Theorem~\ref{gaussSA} theorem below on subbadditivity of
entropy with respect to a Gaussian reference measure is equivalent to
Theorem~\ref{gaussbl}. Hence, it suffices to prove one of the other.

Before stating and proving the subadditivty theorem, we first recall
that for any probability density $\rho$ on $(\R^m,\dd \gamma_m)$, the
entropy of $\rho$, is defined by
$$S(\rho) = -\int_{\R^m}\rho(y)\ln\rho(y) \gamma_m(y)\dd y\ .$$
Note that the relative entropy of $\rho(y)\gamma_m(y)\dd y$ to
$\gamma_m(y)\dd y$ is $-S(\rho)$; in the convention employed here, the
entropy $S$ is concave, and the relative entropy is convex.

\begin{thm}\label{gaussSA} Let $V_1,\dots,V_N$ be $N$ non zero
  subspaces of $\R^n$, and for each $j$, and let $d_j$ denote the
  dimension of $V_j$. Define $P_j$ to be the orthogonal projection of
  $\R^n$ onto $V_j$.  For any probability density $\rho$ on $(\R^n,\dd
  \gamma_n)$, let $\rho_{V_j}$ denote the marginal density on
  $(V_j,\dd \gamma_{d_j})$.  Then, given the numbers $p_j$, $1\le p_j
  < \infty$ for $j = 1,\dots,N$, there exists a finite constant $C$
  such that
\begin{equation}\label{sag}
\sum_{j=1}^N\frac{1}{p_j}S(\rho_{V_j}) \ge S(\rho) - \ln(C)
\end{equation}
holds for all probability densities $\rho$ on 
$(\R^n,\dd \gamma_n)$,  if an only if
\begin{equation}\label{fj}
\sum_{j=1}^N \frac{1}{p_j} P_j \leq I
\end{equation}
and in this case, $\ln(C)=0$.
\end{thm}

We first prove necessity of the condition (\ref{fj}): 

\begin{lm} The condition (\ref{fj}) in Theorem~\ref{gaussSA} is necessary.
\end{lm}

\noindent{\bf Proof:}   It suffices to consider densities of the form 
$$\rho(x) = \exp(b\cdot x - |b|^2/2)\ ,$$
for $b\in \R^n$. 
Then
$$\rho_{V_j}(x) = \exp(P_jb\cdot y - |P_jb|^2/2)\ ,$$
and we compute:
$$S(\rho) = -\frac{|b|^2}{2}\qquad{\rm and}\qquad  S(\rho_{V_j}) = -\frac{|P_jb|^2}{2}\ .$$
Thus
$$\sum_{j=1}^N\frac{1}{p_j}S(\rho_{V_j}) - S(\rho)   = b\cdot \left[  Id_{\R^n} - \sum_{j=1}^N\frac{1}{p_j}P_j\right]b\ ,$$
and evidently this is bounded below if and only if (\ref{fj}) is satisfied. \lanbox

\subsection{Proof of sufficiency}

The sufficiency of the condition (\ref{fj}) will be proved using an
interpolation between an arbitrary density $\rho$ and the uniform
density that is provided by the Mehler semigroup.  (Indeed, Barthe and
Coredero-Erausquin used the Mehler semigroup in their work \cite{BC}
mentioned above, but in a direct proof of the Gaussian B-L inequality
inspired by the heat-flow method introduced in \cite{CLL1}. The heat
flow approach to prove subadditivity inequalities was developed in
\cite{BCM} and \cite{CC}.)

The Mehler semigroup is the strongly continuous semigroup of
positivity preserving contractions on $L^2(\R^n ,\gamma_n(x){\rm d}x)$
whose generator $-{\cal N}$ is given by the Dirichlet form
\begin{equation}\label{gener}
{\cal E}(f,g) = \int_{\R^n} \nabla f^*(x)\cdot \nabla g(x)\gamma_n(x){\rm d}x
\end{equation}
through   $\langle f, {\cal N} g \rangle_{L^2(\gamma_n)} =  {\cal E}(f,g)$,
where $f^*$ is the complex conjugate of $f$.
Integrating by parts, one finds
$${\cal N} = -(\Delta - x\cdot\nabla)\ ,$$
The eigenvalues of ${\cal N}$ are the non-negative integers, and the
eigenfunctions are the Hermite polynomials.  (In certain physical
contexts, the eigenvalues count occupancy of quantum state and ${\cal
  N}$ is called the {\em Boson number operator}.)

There is a simple explicit formula for the $e^{-t{\cal N}}$:
\begin{equation}\label{mfor}
e^{-t{\cal N}}f(x) = \int_{\R^n} f\left(e^{-t}x + \sqrt{1 - e^{2t}}y\right)\gamma_n(y)\dd y\ ,
\end{equation}
which is easily checked. 

Since evidently ${\cal N}1 = 0$, and $e^{-t{\cal N}}$ is self-adjoint,
it also preserves integrals against $\gamma_n(x){\rm d}x$, and hence,
if $\rho$ is any probability density, so is each $\rho_t := e^{-t{\cal
    N}}$.  As one sees from (\ref{mfor}),
\begin{equation}\label{conver}
\lim_{t\to\infty}e^{-t{\cal N}}\rho(x) = 1\ ,
\end{equation}
the uniform probability density on $(\R^n, \gamma_n(x){\rm d}x)$, and
thus the Mehler semigroup provides us with an interpolation between
any probability density $\rho$ and the uniform density $1$.

This interpolation is well-behaved with respect to the operation of
taking marginals: Consider any probability density $\rho$ on
$(\R^n,\gamma_n(x){\rm d}x)$, and any $m$ dimensional subspace $V$ of
$\R^n$. Let $\rho_V$ be the marginal density of $\rho$ as in
Theorem~\ref{gaussSA}. Then of course, we may regard $\rho_V$ as a
probability density on $(\R^n,\gamma_n(x){\rm d}x)$ that is constant
along directions in $V^\perp$.  (Simply compose $\rho_V$ with $P_V$.)
Interpreted this way, so that both $\rho$ and $\rho_V$ are functions
on $\R^N$,
\begin{equation}\label{compo}
\left(e^{-t{\cal N}}\rho\right)_V =  e^{-t{\cal N}} \left(\rho_V\right)\ .
\end{equation}
That is, taking marginals  commutes with the action of the Mehler semigroup. 

The next point to note is that the entropy is monotone increasing
along this interpolation: Differentiating, with $\rho_t = e^{-t{\cal
    N}}\rho$,
$$\frac{{\rm d}}{{\rm d}t} S(\rho_t) =  -\int_{\R^n} \ln(\rho_t)(\Delta  - x\cdot \nabla)\rho_t \gamma_n \dd x = 
\int_{\R^n}\nabla \ln \rho_t \cdot \nabla \rho_t \gamma_n \dd x= {\cal
  E}(\ln \rho_t, \rho_t)\ .$$ For any smooth density $\rho$, $ {\cal
  E}(\ln \rho, \rho) = \int_{\R^n}\nabla \ln \rho \cdot \nabla \rho
\gamma_n \dd x = \int_{\R^n}|\nabla \ln \rho |^2 \rho \gamma_n \dd x
$, and hence $S(\rho_t) $ is strictly increasing for all $t$.
Moreover, since $(x,t) \mapsto |x|^2/t$ is jointly convex on
$\R^n\times \R_+$, $\rho\mapsto {\cal E}(\ln \rho, \rho)$ has a unique
extension as a convex functional the set of all probability densities
on $(\R^n, \gamma_n(x){\rm d}x)$.

\begin{defi}[Entropy Production]\label{eprodg}
The 
{\it entropy production} functional is the convex functional $D(\rho)$ on probability 
densities on $(\R^n, \gamma_n(x){\rm d}x)$ given by
\begin{equation}\label{epfo}
  D(\rho) = \int_{\R^n} \ln \rho(x) {\cal N}\rho(x) \gamma_n(x) {\rm d}x = {\cal E}(\ln \rho, \rho)\ .
\end{equation}
\end{defi}

With this definition,
$$\frac{{\rm d}}{{\rm d}t} S(e^{-t{\cal N}}\rho)  = D(e^{-t{\cal N}}\rho)\ .$$
Now because of (\ref{compo}), for any subspace $V$ of $\R^n$,
$$\frac{{\rm d}}{{\rm d}t} S([e^{-t{\cal N}}\rho]_V)  = D([e^{-t{\cal N}}\rho]_V)\ .$$
Now, since $[e^{-t{\cal N}}\rho]_V$ is constant along directions
orthogonal to $V$, the derivatives in those directions that figure in
$D([e^{-t{\cal N}}\rho]_V)$ are irrelevant; we need only take
derivatives along directions in $V$. This consideration leads to the
definitions of the {\em restricted number operator}, and the {\em
  restricted entropy production}:

Given an $m$ dimensional subspace $V$ of $\R^n$, 
let $P_V$ be the orthogonal projection onto $V$. The restricted number operator ${\cal N}_V$ 
is the self-adjoint operator on 
$L^2(\R^n ,\gamma_n(x){\rm d}x)$ defined through  
\begin{equation}\label{renum}
\langle f, {\cal N}_V g \rangle_{L^2(\gamma_n)} =   \int_{\R^n} \nabla f^*(x)\cdot P_V\nabla g(x)\gamma_n(x){\rm d}x\ ,
\end{equation}
and the  {\em restricted entropy production functional} $D_V(\rho)$ is the convex functional 
given by
\begin{equation}\label{repdef}
D_V(\rho) = \int_{\R^n} \left({\cal N}_V \ln \rho(x) \right)\rho(x) \gamma_n(x) {\rm d}x \ .
\end{equation}

With this definition, $D(\rho_V) = D_V(\rho_V)$, however, there is a
crucial difference between $D_V(\rho)$ and $D(\rho_V)$:

\begin{lm}\label{gaussmon} For any smooth probability density $\rho$ on  
  $(\R^n,\gamma_n(x){\rm d}x)$, and any non $m$ dimensional subspace
  $V$ of $\R^n$, let $\rho_V$ be the corresponding marginal density
  regarded as a probability density on $(\R^n,\gamma_n(x){\rm d}x)$.
  Then
\begin{equation}\label{monpro}
D(\rho_V) \le D_V(\rho)\ .
\end{equation}
\end{lm}

\noindent{\bf Proof:} 
Regard $\rho_{V}$
 as a function on $\R^n$ (by composing it with $P_V$).   
 Assume that $\rho$ is smooth and bounded above and below
 by strictly positive numbers. Notice that since $\rho_V$ is constant 
 constant along directions in $V^\perp$, 
 $${\cal N}\ln \rho_V = {\cal N}_V\ln \rho_V\ ,$$
 and hence

 Then, integrating by parts, and using the definition of $\rho_{V}$
 and the Schwarz inequality, we obtain:
 $$D(\rho_V)=  \int_{\R^n} \left[{\cal N}_V\ln \rho_{V}(x) \right] \rho_{V} (x)
 \gamma_{n}(x) \dd x =  \int_{\R^n} \left[{\cal N}_V\ln \rho_{V} (x)\right] \rho (x)
 \gamma_{n}(x) \dd x\ ,$$
 where we have used the definition of $\rho_V$ to replace the second $\rho_V$ be $\rho$ itself.
 Then, by the definition of ${\cal N}_V$, and the Schwarz inequality,
 
 \begin{eqnarray}\label{sw}
   D(\rho_V) 
   &=&  \int_{\R^n}\left(\nabla \ln \rho_{V} (x)\right)\cdot  P_{V}\nabla  \rho(x)
   \gamma_{n} \dd x\nonumber\\
   &=&  \int_{\R^n}\left(\nabla \ln \rho_{V} (x)\right)\cdot  P_{V}\left(\nabla \ln  \rho(x)\right)
   \rho(x) \gamma_{n}(x) \dd x\nonumber\\
   &\le&  \left(\int_{\R^n}\left|\nabla \ln \rho_{V} (x)\right|^2 \rho(x)
     \gamma_{n}(x) \dd x\right)^{1/2}\left(\int_{\R^n}\left|P_V\nabla
       \ln \rho(x) \right|^2\rho
     \gamma_{n} \dd x\right)^{1/2}\nonumber\\
 \end{eqnarray}
 In the first factor in the last line, we may replace $\rho$ by $\rho_V$ since
 $\left|\nabla \ln \rho_{V} (x)\right|^2$ depends on $x$ only through $P_Vx$. Hence 
 this factor is
 simply $\sqrt{D(\rho_V)}$, and the second factor is $\sqrt{D_V(\rho)}$. 
 \lanbox
 
 The proof we have just given is patterned on the proof of an
 analogous result in the Lebesgue measure case in \cite{CC}, which in
 turn is based on similar arguments in \cite{C} and \cite{BCM}. It is
 somewhat more complicated to adapt the argument to the Clifford
 algebra setting, but this is what we shall do in the next section.
 We are now ready to prove the sufficiency of condition (\ref{fj}):

\begin{lm}\label{gausssuf}The condition (\ref{fj}) in Theorem~\ref{gaussSA} is sufficient.
\end{lm}

\noindent{\bf Proof:} For a probability density $\rho$ on $(\R^n,\dd
\gamma_n)$ $S(\rho)> -\infty$, it is easy to see that
$$\lim_{t\to\infty}S(e^{-t{\cal N}}\rho) = S(1) = 0$$
and hence, $\lim_{t\to\infty}S(e^{-t{\cal N}}(\rho_{V_j})) = 0$ for
each $j=1,\dots,N$.  Therefore, it suffices to show that
$$a(t) := \left[\sum_{j=1}^N\frac{1}{p_j}S(e^{-t{\cal N}} \rho_{V_j}) - S(e^{-t{\cal N}} \rho) \right]$$
is monotone decreasing in $t$.  

Differentiating, and using (\ref{compo}), and then Lemma~\ref{gaussmon},
\begin{eqnarray}
\frac{{\rm d}}{{\rm d}t}a(t) &=&  \left[\sum_{j=1}^N\frac{1}{p_j}D((e^{-t{\cal N}} \rho)_{V_j}) - D(e^{-t{\cal N}} \rho) \right]\nonumber\\
&\le &  \left[\sum_{j=1}^N\frac{1}{p_j}D_{V_j}(e^{-t{\cal N}} \rho) - D(e^{-t{\cal N}} \rho) \right]\nonumber\\
\end{eqnarray}
Now note that by (\ref{repdef}), whenever (\ref{fj}) is satisfied,
$$\sum_{j=1}^N\frac{1}{p_j}D_{V_j}(\sigma) \le  D( \sigma)$$
for any smooth density $\sigma$. Hence the derivative of $\alpha(t)$
is negative for all $t>0$. \lanbox

\section{Generalized subadditivity of the entropy in Clifford algebras}

In this section we shall prove
\begin{thm}\label{cye}  Let $V_1,\dots,V_N$ be $N$ subspaces of $\R^n$, 
  and let $\A_j$ be the Clifford algebra over $V_j$ with the inner
  product $V_j$ inherits from $\R^n$, and let $\A_j$ be equipped with
  its unique tracial state $\tau_j$.  For any probability density
  $\rho\in \A$, let $\rho_{V_j}$ be the induced probability density in
  $\A_j$.  Let $S(\rho) = \tau(\rho\ln \rho)$ and $S(\rho_{V_j}) =
  \tau_j(\rho_{V_j}\ln \rho_{V_j})$

Then 
\begin{equation}\label{clifs1}
\sum_{j=1}^N\frac{1}{p_j} S(\rho_{V_j}) \ge S(\rho)
\end{equation}
for all probability densities $\rho\in \A$ if and only if
\begin{equation}\label{fjc2}
\sum_{j=1}^N\frac{1}{p_j}P_j \le I_{\R^n}\ .
\end{equation}
where 
 $P_j$ is the orthogonal projection onto
$V_j$ in $\R^n$. 
\end{thm}

Granted this result, we have:

\noindent{\bf Proof of Theorem~\ref{cy}:} Theorem~\ref{equiv} and Theorem~\ref{cye}
together prove Theorem~\ref{cy}. \lanbox

We shall now prove
 Theorem~\ref{cye}. 
As before, we begin by proving the necessity of (\ref{fjc2}).

\begin{lm}\label{clnec}
The condition (\ref{fjc2}) in Theorem~\ref{cye} is necessary.
\end{lm}

Proof: For any vector $a = (a_1,\dots,a_n)\in \R^n$, define
$$\rho_a = I + \sum_{j=1}^n a_jQ_j = I + a\cdot Q\ .$$
Then $\rho_a$ is a probability density if and only if $|a| \le 1$.
Indeed, $\rho_a$ has only two eigenvalues, $1\pm |a|$, with equal
multiplicity.

Then $(\rho_a)_{V_j} = I + (P_ja)\cdot Q$, and so $(\rho_a)_{V_j}$ has
only two eigenvalues, $1\pm |P_j a|$, with equal multiplicity.
Therefore,
\begin{equation}\label{ents1}
S(\rho_a) = -\psi(|a|)\qquad{\rm and}\qquad  S((\rho_a)_{V_j}) = -\psi(|P_ja|) \ .
\end{equation}
where $\psi(x)$ is the convex function defined by
\begin{equation}\label{psidef}
\psi(x) =
\begin{cases}
 \frac{1}{2}\left[ (1+x)\ln(1+x) +  (1-x)\ln(1-x)\right] & \text{if $|x| \le 1$.}
\\
\infty & \text{otherwise,}
\end{cases}
\end{equation}

Thus, for (\ref{clifs1}) to hold for each $\rho_a$, $|a|\le 1$, it must be the case that
\begin{equation}\label{ents2}
\sum_{j=1}^N\frac{1}{p_j}\psi(|P_ja|) \le \psi(|a|)\qquad {\rm for\  all\ } a\ {\rm  with} \ |a| \le 1\ .
\end{equation}
Then since $\psi(x) = x^2 + {\cal O}(x^4)$, replacing $a$ by $ta$, $0<
t < 1$, we see that (\ref{fjc2}) must hold. \lanbox

Because of (\ref{ents1}), once we have proved Theorem~\ref{cye}, we
will have a proof of (\ref{ents2}).  However, it is of interest to
have a direct proof of this inequality.

\begin{prop}\label{elem}
The inequality (\ref{ents2}) holds whenever (\ref{fjc2}) is satisfied.
\end{prop}

\noindent{\bf Proof:}  An easy calculation of derivatives shows that 
$$\psi'(x) = {\rm arctanh}(x)\qquad {\rm and}\qquad \psi''(x) = \frac{1}{1-x^2}$$
for $|x| < 1$. 

Now fix any $a$ with $|a| < 1$. Then, for $t>0$, define
$$\eta(t) =   \psi(e^{-t}|a|)  - \sum_{j=1}^N\frac{1}{p_j}\psi(e^{-t}|P_ja|) \ .$$
We have to show that $\eta(t)>0$ for all $t>0$. Since evidently $\lim_{t\to\infty}\eta(t) =0$,
it suffices to show that $\eta'(t) < 0$ for all $t>0$.

Differentiating, we find
$$\eta'(t) = -e^{-t}\left[ |a|\ {\rm arctanh}(e^{-t}|a|) - \sum_{j=1}^N\frac{1}{p_j}
 |P_j a|\ {\rm arctanh}(e^{-t}|P_ja|)\right] := e^{-t}\theta(t)\ .$$

Hence, it suffices to show that $\theta(t)\ge 0$ for all $t>0$. Since once again, 
$\lim_{t\to\infty}\theta(t) =0$,
it suffices to show that $\theta'(t) < 0$ for all $t>0$.
Differentiating, we find
$$\theta'(t) = -e^{-t}\left[\frac{ |a|^2}{1-e^{-2t}|a|^2}- \sum_{j=1}^N \frac{1}{p_j}
  \frac{|P_j a|^2}{1-e^{-2t}|P_ja|^2}\right] \ .$$ Multiplying through
by $e^{-t}$, and absorbing a factor of $e^{-t}$ into $a$, it suffices
to show that
\begin{equation}\label{rat1}
\frac{ |a|^2}{1-|a|^2}\ge  \sum_{j=1}^N \frac{1}{p_j}
\frac{|P_j a|^2}{1-|P_ja|^2} 
\end{equation}
for all $|a|\le 1$.  However, since $|a| \ge |P_ja|$,
$$\frac{|P_j a|^2}{1-|a|^2} \ge \frac{|P_j a|^2}{1-|P_ja|^2}\ ,$$
and thus (\ref{rat1}) follows from (\ref{fjc2}). \lanbox

We are now in a position to give an elementary proof of
Theorem~\ref{cy} in the special case that each $V_j$ is one
dimensional. As explained in Example~\ref{cliff}, it suffices in this
case to prove the following:

\begin{prop}\label{onedimc} Suppose $\{u_1,\dots,u_N\}$ is any set of $N$ unit vectors in $\R^n$, and 
$\{p_1,\dots,p_N\}$ is any set of $N$ positive numbers such that
\begin{equation}\label{fjfj}
\sum_{j=1}^N c_j u_j \otimes u_j =I_{\R^n}\ .
\end{equation}
Then for any $b = (b_1,\dots,b_N)\ in \R^N$, 
\begin{equation}\label{cosh2b}
\ln\cosh\left(\left|{\textstyle \sum_{j=1}^N b_ju_j}\right|\right) \le \sum_{j=1}^N\frac{1}{p_j}
\ln\cosh(p_jb_j) \ .
\end{equation}
\end{prop}

\noindent{\bf Proof:} Let $\psi^*(x)$ denote the function $\psi^*(x) =
\ln \cosh(x)$, $x\in \R$. The notation is meant to indicate the well
known fact, easily checked, that $\psi^*$ is the Legendre transform of
the function $\psi$ defined in (\ref{psidef}).

Now, given a set of $N$ orthogonal projections $\{P_1,\dots,P_N\}$
satisfying (\ref{fjc2}), we may make any choice of a unit vector $u_j$
from the range of $P_j$, and then the $N$ unit vectors
$\{u_1,\dots,u_N\}$ will satisfy (\ref{fjfj}). Conversely, given any
set of $N$ unit vectors $\{u_1,\dots,u_N\}$ that satisfy (\ref{fjfj}),
we may take $P_j = u_j\otimes u_j$, and then (\ref{fjc2}) is
satisfied.  Hence, we suppose we are given a a set of $N$ orthogonal
projections $\{P_1,\dots,P_N\}$ satisfying (\ref{fjc2}), and for each
$j$, $u_j$ is a unit vector in the range of $P_j$.

Then for any $b \in \R^n$,
\begin{eqnarray}
\psi^*\left(\left|{\textstyle \sum_{j=1}^N b_ju_j}\right|\right) 
&=& \sup_{a\in \R^n}\left\{ a \cdot  \sum_{j=1}^N b_ju_j  - \psi(|a|)\right\}\nonumber\\
&=& \sup_{|a| \le 1}\left\{ \sum_{j=1}^N P_ja \cdot   b_ju_j  - \psi(|a|)\right\}\nonumber\\
&\le& \sup_{|a| \le 1}\left\{ \sum_{j=1}^N P_ja \cdot   b_ju_j  - \sum_{j=1}^N\frac{1}{p_j}\psi(|P_ja|)\right\}\nonumber\\
&\le& \sup_{|a| \le 1}\left\{ \sum_{j=1}^N |P_ja|   | b_j|  - \sum_{j=1}^N\frac{1}{p_j}\psi(|P_ja|)\right\}\nonumber\\
&=& \sup_{|a| \le 1}\left\{ \sum_{j=1}^N\frac{1}{p_j} \bigl[ |P_ja|   p_j| b_j|  -\psi(|P_ja|)\bigr]\right\}\nonumber\\
\end{eqnarray}
where the first inequality is from (\ref{ents2}), and the second is from Schwarz. Then, by the definition of the Legendre transform,  for any $a$,
$$\psi^*(p_jb_j) \ge  |P_ja|  ( p_j| b_j|)  -\psi(|P_ja|)\ ,$$
we obtain
$$\psi^*\left(\left|{\textstyle \sum_{j=1}^N b_ju_j}\right|\right)   \le \sum_{j=1}^N\frac{1}{p_j}
\psi^*(p_jb_j)\ ,$$
which is (\ref{cosh2b}). \lanbox

We now prove Theorem~\ref{cye} in full generality.  This gives another
proof of the last two propositions, but by less elementary means.  The
proof will follow the basic pattern of the proof of
Theorem~\ref{gaussSA}, and use the Clifford algebra analog of the
Mehler semigroup.  This is the so-called Clifford--Mehler semigroup,
about which we now recall a few relevant facts.
 
\subsection{About the Clifford--Mehler semigroup}

There is also a differential calculus in the Clifford algebra. Let
$Q_1,\dots,Q_n$ be any set of $n$ generators for the Clifford algebra
$\CC$ over $\R^n$.  For $A\in \CC$, define
$$\nabla_i(A) = \frac{1}{2}\left[Q_iA - \Gamma(A)Q_i\right]\ ,$$
where $\Gamma$ is the {\em  grading operator} on $\CC$: That is, using the notation in (\ref{multii}),
$$\Gamma(Q^{\alpha}) = (-1)^{|\alpha|}Q^\alpha\ .$$
One computes that $\nabla_i(Q^\alpha) = 0$ is $\alpha(i) = 0$, and
otherwise, $\nabla_i(Q^\alpha) = 0$ is what one gets by anti-commuting
the factor of $Q_i$ through to the left, and then removing it. In this
sense it is like a differentiation operator, and what is more, it is a
skew derivation on $\CC$, which means that for all and $A$ and $B$ in
$\CC$, $\nabla_j(AB) = \nabla_j(A) B + \Gamma(A)\nabla_j(B)$.

The Clifford algebra analog of the Gaussian energy integral
(\ref{gener}) is given by
\begin{equation}\label{generc}
{\cal E}(A,B) =\tau\left( \sum_{j=1}^n \nabla_j A^*\nabla_j B\right)\ ,
\end{equation}
for all $A,B\in \CC$. This is the {\em Clifford Dirichlet form} studied by Gross.
Then, the Fermionic number operator, also denoted ${\cal N}$, is defined by
$$  {\cal E}(A,B) = \tau(A^* {\cal N}(B))\ .$$
It is easy to see that the spectrum of ${\cal N}$ consists of the non
negative integers $\{0,1,\dots,n\}$ and that
\begin{equation}\label{deg}
{\cal N}Q^\alpha = |\alpha|Q^\alpha\ .
\end{equation}

The Clifford Mehler semigroup is then given by $e^{-t{\cal N}}$.  It
is clear from this definition, (\ref{taudef}) and (\ref{deg}) that for
any $A\in \CC$, $\lim_{t\to\infty}e^{-t{\cal N}}(A) = \tau(A)I$.  Thus
for any probability density $\rho$ in $\CC$,
$$t\mapsto \rho_t =  e^{-t{\cal N}}(\rho)$$
provides an interpolation between $\rho$ and $I$, and each $\rho_t$ is
a probability density.  This corresponds exactly to the Mehler
semigroup interpolation that was used to prove Theorem~\ref{gaussSA},
and we shall use it here in the same way, though some additional
complications shall arise.

${\cal N}$ does not depend on the choice of the set of generators
$Q_1,\dots,Q_n$. Indeed, if $\{u_1,\dots,u_n\}$ is any orthonormal
basis of $\R^n$, and we define $\widetilde Q_j = u_j\cdot Q\qquad j
=1,\dots,n$, then the Clifford Dirichlet form that one obtains using
this basis to define the derivatives is the same as the original.

In particular, given an $m$ dimensional subspace $V$ of $\R^n$, we may
choose $\{u_1,\dots,u_n\}$ so that $\{u_1,\dots,u_m\}$ is an
orthonormal basis for $V$, and then the first $m$ generators will be a
set of generators for $\CC_V$.  We then define the {\em reduced
  Clifford Dirichlet form} ${\cal E}_V$ by
\begin{equation}\label{generc2}
{\cal E}_V(A,B) =\tau\left( \sum_{i,j=1}^n \nabla_i A^*[P_V]_{i,j} \nabla_j B\right)\ ,
\end{equation}
where $[P_V]_{i,j}$ is the $i,j$th entry of the $n\times n$ matrix for $P_V$.
The restricted number operator ${\cal N}_V$ is then the self-adjoint operator on $L^2(\CC)$
given by
$\tau(A^* {\cal N}_V(B)) =  {\cal E}_V(A,B)$.

Now, for any probability density $\rho$ in $\CC$ let $\rho_V$ be the
corresponding marginal density regarded as an operator in $\CC$ by
identifying it with $\phi_V(\rho_V)$, where $\phi_V$ is the canonical
embedding of $\CC(V)$ into $\CC(\R^n)$. Then it is an easy consequence
of the definitions that
\begin{equation}\label{clcond}
\left(e^{-t{\cal N}}\rho\right)_V =    e^{-t{\cal N}} \left(\rho_V\right) = e^{-t{\cal N}_V} \left(\rho_V\right)\ .
\end{equation}

Also, under the condition (\ref{fjc2}), it is easy to see that 
\begin{equation}\label{fjc3}
\sum_{j=1}^N\frac{1}{p_j}{\cal N}_{V_j}  \le {\cal N}\ .
\end{equation}

Finally, we introduce {\it entropy production} $D(\rho)$:  With $\rho_t := e^{-t{\cal N}}\rho$
,we differentiate and  find
$$\frac{{\rm d}}{{\rm d}t} S\left(\rho_t\right) = 
\tau\left(\ln(\rho_t){\cal N}(\rho_t)\right) = {\cal E}(\ln(\rho_t),\rho_t)\ .$$

\begin{defi}[Entropy Production] The {\it entropy production}
  functional at a probability density $\rho$ is the functional defined
  by
\begin{equation}\label{epformb}
D(\rho) = \tau\left(\ln(\rho){\cal N}(\rho)\right) = {\cal E}(\ln(\rho),\rho)\ .
\end{equation}

Given an $m$ dimensional subspace $V$ of $\R^n$, the {\it restricted
  entropy production} functional at a probability density $\rho$ is
the functional defined by
\begin{equation}\label{epform}
D_V(\rho) = \tau\left(\ln(\rho){\cal N}_V(\rho)\right) = {\cal E}_V(\ln(\rho),\rho)\ .
\end{equation}
\end{defi}

The following lemma is the basis of our proof of the sufficiency of
(\ref{fjc2}).  In the course of proving it, we shall see that both
$D(\rho)$ and $D_V(\rho)$ are convex functionals, which is somewhat
less obvious than in the Gaussian case.

\begin{lm}\label{clmon} For any any probability density $\rho$  in $\CC(\R^n)$, and any
  $m$ dimensional subspace $V$ of $\R^n$, let $\rho_V$ be the
  corresponding marginal probability density regarded as a probability
  density in $\CC(\R^n)$. Then
$$D(\rho_V) \le D_V(\rho)\ .$$
\end{lm}

\noindent{\bf Proof:} We choose an orthonormal basis
$\{u_1,\dots,u_n\}$ for $\R^n$ such that $\{u_1,\dots,u_m\}$ is an
orthonormal basis for $V$. Without loss of generality, we may suppose
that $\{u_1,\dots,u_n\}$ is the standard basis so that
$\{Q_1,\dots,Q_m\}$ is a set of generators for $\CC(V)$.  Then,
\begin{equation}\label{generc3}
{\cal E}(A,B) =\tau\left( \sum_{j=1}^n \nabla_j A^*\nabla_j B\right)\qquad{\rm and}\qquad
{\cal E}_V(A,B) =\tau\left( \sum_{j=1}^m \nabla_j A^*\nabla_j B\right)\ .
\end{equation}
It will be convenient to define
${\cal N}_j =  \nabla_j^*\nabla_j\qquad j = 1,\dots,n$.
Then we have
\begin{equation}\label{gen44}
{\cal N} = \sum_{j=1}^n{\cal N}_j   \qquad{\rm and}\qquad
{\cal N}_V = \sum_{j=1}^m {\cal N}_j  \ ,
\end{equation}
and so
\begin{equation}\label{pro3}
D_V(\rho) = \sum_{j=1}^m \tau\left( \ln \rho, {\cal N}_j \rho\right)\ .
\end{equation}
Since 
${\displaystyle {\cal N}_j Q^{\alpha} = 
\begin{cases}
Q^{\alpha} & \text{if $\alpha(j) =1$,}
\\
0 & \text{$\alpha(j)=0$,}
\end{cases}
}$, each ${\cal N}_j$ is an orthogonal projection, and so (\ref{pro3})
can be rewritten as
\begin{equation}\label{pro3b}
D_V(\rho) = \sum_{j=1}^m \tau\left( {\cal N}_j(\ln \rho), {\cal N}_j \rho\right)\ .
\end{equation}

To proceed, we use a formula of Gross \cite{G75} for ${\cal N}_j f(A)$
where $A\in \CC(\R^n)$, and $f$ is a continuous function. To write
down Gross's formula, first write $A = B + Q_jC$ where both $B$ and
$C$ are linear combinations of the $Q^\alpha$ with $\alpha(j) =0$.
Then define $\widehat A = B - Q_jC$.  Notice that if $\rho$ is a
probability density, then $\widehat\rho$ is again a probability
density.  Gross's formula is
$${\cal N}_j f(A) = \frac{1}{2}\left[ f(A) - f(\widehat A\ )\right]\ .$$

To prove this formula, notice that there is a unitary operator $U$
such that $\widehat A = UAU^*$.  (If the dimension $n$ is odd, one can
take $U$ to be the product, in some order, of all of the $Q_k$ for
$k\ne j$; if the dimension is even, one can add another generator.)
Therefore,
$$\widehat{f(A)} = U f(A) U^* = f(UAU^*) = f(\widehat A\ )\ .$$
Using this together with the fact that for any $A\in \A$, ${\cal N}_jA = (1/2)[A - \widehat{A}]$,
we obtain Gross's formula, which we now apply as follows:

\begin{eqnarray}
\tau\left( {\cal N}_j(\ln \rho) {\cal N}_j \rho\right) &=&
\frac{1}{4}\tau\left( \left[ \ln(\rho) - \ln(\widehat\rho)\right]  \left[ \rho - \widehat\rho\right]  \right)\nonumber\\
&=&
\frac{1}{4}\tau\left(  \ln(\rho)\left[ \rho - \widehat\rho\right] \right)  +  \frac{1}{4}\tau\left(  \ln(\widehat\rho)\left[ \widehat\rho - \rho\right] \right)\nonumber\\
&=& \frac{1}{4} H[\rho|\widehat\rho\ ] + \frac{1}{4} H[\widehat\rho\ |\rho] \nonumber\\
\end{eqnarray}
where $H[\rho|\sigma] = \tau \rho(\ ln \rho - \ln \sigma)$ is the {\em
  relative entropy} of $\rho$ with respect to $\sigma$.  As is well
known, $(\rho,\sigma) \mapsto H(\rho|\sigma)$ is jointly convex, and
hence
$$\rho \mapsto  \tau\left( (\ln \rho){\cal N}_j \rho\right) $$
is convex. Furthermore, by the fundamental monotonicity property of
the relative entropy under conditional expectations \cite{U2},
$$H(\rho_V|\sigma_V) \le H(\rho|\sigma)$$
for any two probability densities $\rho$ and $\sigma$. It follows that
$ \tau\left( (\ln \rho_V) {\cal N}_j \rho_V\right) \le \tau\left( (\ln
  \rho) {\cal N}_j \rho\right)$.  Summing on $j$ from $1$ to $m$, we
find
$$D(\rho_V) = D_V(\rho_V) = \sum_{j=1}^m\tau\left( (\ln \rho_V) {\cal N}_j \rho_V\right) \le   
\sum_{j=1}^m \tau\left( (\ln \rho) {\cal N}_j \rho\right) = D_V(\rho)\
.$$ \lanbox

\subsection{Proof of the sufficiency}    

\begin{lm}\label{clsuf}The condition (\ref{fj}) in Theorem~\ref{gaussSA} is sufficient.
\end{lm}

\noindent{\bf Proof:}  For   a probability density $\rho$ in  $\CC(\R^n)$
it is easy to see that 
$$\lim_{t\to\infty}S(e^{-t{\cal N}}\rho) = S(1) = 0$$
and hence,
$\lim_{t\to\infty}S(e^{-t{\cal N}}(\rho_{V_j}))  = 0$ for each $j=1,\dots,N$.
Therefore,
it suffices to show that 
$$a(t) := \left[\sum_{j=1}^N\frac{1}{p_j}S(e^{-t{\cal N}} \rho_{V_j}) - S(e^{-t{\cal N}} \rho) \right]$$
is monotone decreasing in $t$.  

Differentiating, and using (\ref{clcond}), and then Lemma~\ref{clmon},
\begin{eqnarray}
\frac{{\rm d}}{{\rm d}t}a(t) &=&  \left[\sum_{j=1}^N\frac{1}{p_j}D((e^{-t{\cal N}} \rho)_{V_j}) - D(e^{-t{\cal N}} \rho) \right]\nonumber\\
&\le &  \left[\sum_{j=1}^N\frac{1}{p_j}D_{V_j}(e^{-t{\cal N}} \rho) - D(e^{-t{\cal N}} \rho) \right]\nonumber\\
\end{eqnarray}
Now note that by (\ref{repdef}), whenever (\ref{fj}) is satisfied,
${\displaystyle \sum_{j=1}^N\frac{1}{p_j}D_{V_j}(\sigma) \le D(
  \sigma)}$ for any smooth density $\sigma$. Hence the derivative of
$\alpha(t)$ is negative for all $t>0$. \lanbox

Notice that the proof is almost identical, symbol for symbol, with
that of the corresponding proof in the Gaussian case. The main
difference of course is that the proof of the main lemma,
Lemma~\ref{clmon}, is considerably more intricate than that of its
Gaussian counterpart.

\noindent{\bf Proof of  Theorem~\ref{cye}:}  This now follows immediately from Lemma~\ref{clnec} and \ref{clsuf}. \lanbox


\bigskip

\begin{thebibliography}{99}

\footnotesize{

\bibitem{B} 
[4] F. Barthe,  \textit{ On a reverse form of the Brascamp-Lieb inequality}, Invent. Math. {\bf 134} (1998), no. 2, pp. 335--361.

\bibitem{BC}   F.~Barthe and D.~Cordero--Erausquin
 \textit{ Inverse Brascamp-Lieb inequalities along the Heat equation}, in {\em Geometric Aspects of Functional Analysis} (2002-2003), LNM {\bf 1850}, Springer, 2004, pp. 65-71.

\bibitem{BCM} F.~Barthe, D.~Cordero--Erausquin and B. Maurey:  \textit{ Entropy of spherical marginals and related inequalities}, J. Math. Pures Appl.  {\bf 86},  no. 2, 2006 pp. 89--99.


\bibitem{Beckner} W.~Beckner, {\it Inequalities in Fourier
series}, Ann. of Math.  {\bf 102},  1975, pp. 159-182.

\bibitem{BL} H.J.~Brascamp and E.H.~Lieb:  \textit{The best constant in Young's inequality and its generalization to more than three functions}, Adv. in Math. {\bf 20},  1976,
pp. 151-173




\bibitem{BW} R.~Brauer, R., H~Weyl:  \textit{Spinors in $n$ dimensions}, Am. Jour.
Math., {\bf 57}, 1935 pp. 425-449.



\bibitem{BCCT1}
J.~Bennett, A.~Carbery, M.~Christ and T.~Tao:   \textit{ The Brascamp-Lieb inequalities: finiteness, structure, and extremals}, Preprint.

\bibitem{BCCT2}
J.~Bennett, A.~Carbery, M.~Christ and T.~Tao:
 \textit{ Finite bounds for Holder-Brascamp-Lieb multilinear inequalities}, Preprint.

\bibitem{C}E.A.~Carlen:  \textit{Superadditivity of Fisher information and logarithmic Sobolev 
inequalities}, Jour. of Func. Analysis,  {\bf 101}, 1991, pp. 194--211.

\bibitem{CC}E.A. Carlen and D Cordero--Erausquin: \textit{
Subadditivity of the entropy and its relation to Brascamp-Lieb type inequalities},  to appear in Geom. Funct. Anal. (2008).

\bibitem{CL1} E.A.~Carlen and E.H.~Lieb:  \textit{ Optimal Hypercontractivity for Fermi Fields and Related  Non-commutative Integration Inequalities}, Comm. Math. Phys.,  {\bf 155}, 1993 pp. 27-46.

\bibitem{CL2} E.A.~Carlen and E.H.~Lieb:   \textit{ Optimal two-uniform convexity and fermion 
hypercontractivity}, Proceedings of the Kyoto Conference in Honor of Araki's 
60th birthday, Kluwer Academic Publishers, Amsterdam, 1993 

\bibitem{CLL1} E.A.~Carlen,  E.H.~Lieb and M.~Loss:  \textit{ A sharp form of Young's inequality on $S^N$ and related entropy inequalities},
Jour. Geom. Analysis {\bf 14}, 2004,  pp. 487-520 

\bibitem{CLL2} E.A.~Carlen,  E.H.~Lieb and M.~Loss:  \textit{ A inequality of Hadamard type for permanents}, Meth. and Appl. of Analysis,
{\bf 13}, no. 1, 2006 pp. 1-17 

\bibitem{Di53} 
J.~Dixmier:  \textit{Formes lin\'eaires sur un anneau d'op\'erateurs},
Bull. Soc. Math. France {\bf 81}, 1953, pp.  222--245

\bibitem{FK} T.~Fak and H.~Kosaki: \textit{Generalized $s$-numbers of $\tau$--measurable operators.} Pacific J. Math.,  {\bf 123}, 1986, pp. 269-300, 

\bibitem{Go} S.~Golden:  \textit{Lower bounds for Helmholtz functions}. Phys. Rev. {\bf 137B},  1965, pp. 1127-1128 

\bibitem{G75} L.~Gross:  \textit{Hypercontractivity and logarithmic Sobolev
inequalities for the Clifford-Dirichlet form}, Duke Math. J., {\bf 43}, 1975, pp.
 383-396 

\bibitem{H} U.~Haagerup:  \textit{$L^p$ spaces associated with an arbitrary von Neumann algebra}, Alg\`ebras d'ope\^erateurs et leurs applications en physique math\^ematique (Colloque CNRS, Marseille, juin 1977) Editions du CNRS, Paris, 1979,  pp. 383-396.

\bibitem{HLP} G.H.~Hardy, J.E. Littlewood and G. P\'olya:  \textit{Inequalities}. Cambridge University Press, Cambridge, 1934

\bibitem{K2} H.~Kosaki: \textit{Applications of the complex
interpolation method to a von Neumann algebra}: noncommutative $L^p$
-spaces.} J. Funct. Anal., {\bf 56}, 1984,  no. 1, pp. 29-78.

\bibitem{L75} E.H.~Lieb: \textit {Some convexity and subadditivity
properties of entropy}, Bull. Amer. Math. Soc. {\bf 81}, 1975, pp. 1-13.

\bibitem{L90} E.H.~Lieb: \textit {Gaussian kernels have only Gaussian maximizers}, Invent. Math. {\bf 102}, 1980, pp. 179-208 

\bibitem{LR}  E.H.~Lieb and M.B.~Ruskai: \textit{Proof of the strong
subadditivity of quantum-mechanical entropy}, J. Math. Phys. {\bf 14}, 1973, pp.
1938-1941.


\bibitem{N74} E.~Nelson:  \textit{Notes on non-commutative integration}, Jour. Funct.
Analysis, {\bf 15}, 1974, pp. 103-116 


\bibitem{P88} D.~Petz:  \textit{A variational expression for the relative entropy}, Comm. Math.
Phys., {\bf 114}, 1988, pp.   345--349 

\bibitem{P1}G.~Pisier. \textit{Non-commutative vector valued $L^p$--spaces
and completely p--summing maps} Ast\'erisque, {\bf 247}, Math. Soc. of
France, Paris, 1998.


\bibitem{R} T.~Rockafellar: \textit{Conjugate duality and optimization},
Vol. 16,  \textit{Regional conference series in applied mathematics},
SIAM, Philadelphia, 1974

\bibitem{RU1} M.B.~Ruskai, \textit{Inequalities for quantum entropy: A
review with conditions for equality}, J. Math. Phys. {\bf 43}, 2005, pp. 4358-4375
(2002). Erratum {\it ibid} {\bf 46}, pp. 019901 

\bibitem{S53}  I.E.~Segal: \textit{A non-commutative extension of abstract
integration}, Annals of Math., {\bf 57}, 1953, pp.  401--457 

\bibitem{S56}   I.E.~Segal.: \textit{Tensor algebras over Hilbert spaces II}, 
Annals of Math., {\bf 63}, 1956, pp.  160--175 

\bibitem{S65}   I.E.~Segal.: \textit{Algebraic integration theory}, 
Bull. Am. Math. Soc.., {\bf 71}, 1965,  no. 3,  pp. 419-489 

\bibitem{St59}  W.~Stinespring: \textit{Integration theorems for gages and duality for unimodular groups}, 
Trans. Am. Math. Soc., ., {\bf 90}, 1958,  pp. 15-26 

\bibitem{T}  C.~Thompson: \textit{Inequality with application in statistical mechanics} J. Math. Phys. 
{\bf 6}, 1812--1813 (1965)

\bibitem{U2} A.~Uhlmann: \textit{Relative entropy and the
Wigner-Yanase-Dyson-Lieb concavity in an interpolation theory}
Commun. Math. Phys., {\bf 54}, 1977, pp.  21-32 

\bibitem{Um1} H.~Umegaki: \textit{Conditional expectation in operator algebras I},
Tohoku Math. J., {\bf 6}, 1954, pp.  177-181   

\bibitem{VDN} D. Voiculescu, K.~Dykema and A.~Nica: \textit{Free random variables},
CRM Monograph Series,  {\bf 1}, Am. Math. Soc.,  Providence, R.I., 1992)

\bibitem{Y12} W.H..~Young: \textit{On the multiplication of successions of Fourier constants},
Proc. Royal soc. A., {\bf 97}, 1912,  pp. 331-339  

\bibitem{Y12b} W.H..~Young: \textit{Sur la g\'en\'eraliation du th\'eore\`me du Parseval},
\textit{Comptes rendus.}, {\bf 155}, 1912, pp. 30-33  

\end{thebibliography}
\end{document}